\documentclass[graybox]{svmult}

\usepackage[T1]{fontenc}

\usepackage{type1cm}        
\usepackage{makeidx}         
\usepackage{graphicx}        
\usepackage{multicol}        
\usepackage[bottom]{footmisc}

\usepackage{newtxtext}       %
\usepackage[varvw]{newtxmath}       

\makeindex             
                       
\usepackage{hyperref}
\usepackage{array}
\newcolumntype{K}[1]{>{\centering\arraybackslash$}p{#1}<{$}}

\usepackage[dvipsnames,svgnames,table]{xcolor}

\begin{document}

\title*{Theoretical properties of the eigenvector method}
\author{S\'andor Boz\'oki\orcidID{0000-0003-4170-4613} and\\ L\'aszl\'o Csat\'o\orcidID{0000-0001-8705-5036}}
\institute{S\'andor Boz\'oki \at HUN-REN Institute for Computer Science and Control (SZTAKI), Hungary, 1111 Budapest, Kende street 13--17., \email{bozoki.sandor@sztaki.hun-ren.hu} \at Corvinus University of Budapest (BCE), Hungary, 1093 Budapest, F{\H o}v\'am square 8. \and
L\'aszl\'o Csat\'o \at HUN-REN Institute for Computer Science and Control (SZTAKI), Hungary, 1111 Budapest, Kende street 13--17., \email{laszlo.csato@sztaki.hun-ren.hu} \at Corvinus University of Budapest (BCE), Hungary, 1093 Budapest, F{\H o}v\'am square 8.}
%
%
\maketitle


\abstract{A classical proposal to derive weights from a pairwise comparison matrix is the right eigenvector. The literature has identified some potential weaknesses of this method in previous decades. This chapter discusses five of these issues. First, right-left asymmetry emerges because of the difference between the right and inverse left eigenvectors. Second, group incoherence for choice means that, in group decision-making problems, the ranking given by the aggregated individual weight vectors is not guaranteed to coincide with the ranking derived from the aggregated pairwise comparison matrix. Third, the ranking based on the right eigenvector may depend on the intensity of the preferences, represented by taking a positive power of all comparisons. Fourth, both the ranking position and the normalised weight of an object might change counter-intuitively after modifying a particular comparison. Fifth, the right eigenvector is not necessarily Pareto efficient: a dominating weight vector that approximates each pairwise comparison at least as well, with an improvement in at least one position, could exist. All violations of the theoretical properties are highlighted by illustrative examples. We also present several open questions in order to inspire future research.}



\section{Introduction} \label{Sec1}

The eigenvector method, suggested to derive weights from a pairwise comparison matrix by the founder of the Analytic Hierarchy Process (AHP) methodology \cite{Saaty1977, Saaty1980}, remains the most popular weighting procedure. However, it does not satisfy some attractive theoretical properties. This chapter presents five issues that have been identified in the literature: right-left asymmetry (Section~\ref{Sec3}), group incoherence for choice (Section~\ref{Sec4}), sensitivity to the intensity of preferences (Section~\ref{Sec5}), non-monotonicity (Section~\ref{Sec6}), and Pareto inefficiency (Section~\ref{Sec7}).
Even though we do not discuss all details of these potential weaknesses, the interested reader is provided with informative illustrations and further references in each case.

\section{Preliminaries} \label{Sec2}

Denote a set of alternatives by integers from 1 to $n$.
The $n \times n$ positive matrix $\mathbf{A} = \left[ a_{ij} \right]$ is called (multiplicative) \emph{pairwise comparison matrix} if $a_{ji} = 1/a_{ij}$ for all $1 \leq i,j \leq n$, which naturally implies $a_{ii} = 1$ for all $1 \leq i \leq n$.
Entry $a_{ij}$ measures the degree to which alternative $i$ is preferred to alternative $j$.

Pairwise comparisons are often used to derive the weights of the alternatives, collected in a positive vector $\mathbf{w} = \left[ w_{i} \right]$.
The ratio of weights $w_i$ and $w_j$ is expected to approximate $a_{ij}$ as well as possible. Thus, the weight vector is defined only up to multiplication, and is usually normalised, for example, by $\sum_{i=1}^n w_{i} = 1$.
A \emph{weighting method} associates a weight vector with any pairwise comparison matrix.

A pairwise comparison matrix is said to be \emph{consistent} if the condition $a_{ik} = a_{ij} a_{jk}$ holds for all $1 \leq i,j,k \leq n$; otherwise, it is said to be \emph{inconsistent}.
A pairwise comparison matrix $\mathbf{A} = \left[ a_{ij} \right]$ is consistent if and only if there exists a weight vector $\mathbf{w} = \left[ w_{i} \right]$ such that $a_{ij} = w_i / w_j$ for all $1 \leq i,j \leq n$.

Any reasonable weighting method should give this unique weight vector for a consistent pairwise comparison matrix; this axiom is called correct result in the consistent case by \cite{Fichtner1984} and correctness by \cite{Csato2019a}.
On the other hand, more than one weight vector would make sense in the inconsistent case. Indeed, many weighting methods have been suggested; see \cite{ChooWedley2004} for a thorough discussion of them.

The most widely used procedures are the \emph{(row) geometric mean (logarithmic least squares) method} \cite{WilliamsCrawford1980, CrawfordWilliams1985, DeGraan1980, deJong1984, Rabinowitz1976}, and the (row) \emph{eigenvector method} \cite{Saaty1977, Saaty1980}.
Even though the geometric mean is supported by strong axiomatic results \cite{Fichtner1984, Barzilai1997, LundySirajGreco2017, Csato2018b, Csato2019a}, the Analytic Hierarchy Process (AHP) mainly uses the eigenvector method since it was developed and recommended by \emph{Thomas L.~Saaty}.

The eigenvector method associates the right eigenvector $\mathbf{w}^{(\mathit{EM})} (\mathbf{A})$ as the weight vector with any pairwise comparison matrix $\mathbf{A}$:
\[
\mathbf{A} \mathbf{w}^{(\mathit{EM})}(\mathbf{A}) = \lambda_{\max}(\mathbf{A}) \mathbf{w}^{(\mathit{EM})}(\mathbf{A}),
\]
where $\lambda_{\max}(\mathbf{A})$ is the dominant/maximal/principal/Perron eigenvalue of matrix $\mathbf{A}$.

The (row) geometric mean, or logarithmic least squares, method associates the vector given by the geometric mean of row elements with a pairwise comparison matrix $\mathbf{A}$:
\[
w_i^{(\mathit{GM})}(\mathbf{A}) = \frac{\prod_{j=1}^n a_{ij}^{1/n}}{\sum_{k=1}^n \prod_{j=1}^n a_{kj}^{1/n}}.
\]

Since the priority vectors of different weighting methods may differ for inconsistent matrices, it is important to quantify the level of inconsistency. A number of indices have been suggested for this purpose; \cite{Brunelli2018} provides an overview of them.
\cite{Saaty1977} proposed the first inconsistency index---\cite{Saaty1977} uses the term consistency index but it is better called inconsistency index since consistency is a binary concept---together with the eigenvector method: for any pairwise comparison matrix $\mathbf{A}$, the \emph{inconsistency index} is
\[
\mathit{CI} (\mathbf{A}) = \frac{\lambda_{\max}(\mathbf{A})-n}{n-1}.
\]
$\mathit{CI} (\mathbf{A}) \geq 0$, and the equality holds if and only if the pairwise comparison matrix is consistent.
The \emph{random index} $\mathit{RI}_n$ is the average $\mathit{CI}$ of a large number of pairwise comparison matrices, whose entries are randomly and uniformly generated from the so-called Saaty scale:
\begin{equation} \label{Saaty_scale}
\left\{ 1/9,\, 1/8,\, \dots ,\, 1/2,\, 1,\, 2,\, \dots ,\, 9 \right\}.
\end{equation}
These values are reported in several studies, see, for example, \cite{BozokiRapcsak2008}.
Finally, the ratio $\mathit{CI} / \mathit{RI}_n$ is called the \emph{inconsistency ratio} $\mathit{CR}$. According to \cite{Saaty1977}, the level of inconsistency can be tolerated if the value of $\mathit{CR}$ remains below the threshold $0.1$.

\section{Right-left asymmetry} \label{Sec3}

In this section, the right eigenvector is denoted by $\mathbf{w}^{(R)} = \left[ w_i^{(R)} \right]$.

The dominant eigenvalue $\lambda_{\max}(\mathbf{A})$ has an associated left eigenvector $\mathbf{w}^{(L)} = \left[ w_i^{(L)} \right]$, too, satisfying the system of linear equations $\mathbf{w}^{(L)} \mathbf{A} = \lambda_{\max}(\mathbf{A}) \mathbf{w}^{(L)}$. In the consistent case, there exists a constant $c$ such that $w_i^{(L)} = c / w_i^{(R)}$ for all $1 \leq i \leq n$, and the entry-wise inverse of the left eigenvector, $\mathbf{w}^{(-L)}$, is the right eigenvector.

Therefore, the left eigenvector can also be used to derive the priorities if one takes the entry-wise inverse $\mathbf{w}^{(-L)}$ of it.
The method based on the right eigenvector is not ``better'' than the method based on the (inverse) left eigenvector in any sense \cite{JohnsonBeineWang1979}. The entry $a_{ij}$ of a pairwise comparison matrix $\mathbf{A}$ answers the question \emph{``How much better/more important is alternative $i$ than alternative $j$?''}---but the reciprocal question of \emph{``How much is alternative $j$ worse/less important than alternative $i$?''} can be equivalently asked. This leads to the transposed matrix $\mathbf{A}^\top$, and the right eigenvector of $\mathbf{A}^\top$ coincides with the entry-wise inverse of the left eigenvector of $\mathbf{A}$ \cite{BozokiRapcsak2008}.

The natural question is whether the two approaches lead to the same weights or ranking.
If the number of alternatives is $n=3$, the row geometric mean weight vector $\mathbf{w}^{(\mathit{GM})}$ coincides with the right eigenvector $\mathbf{w}^{(R)}$ \cite[p.~393]{CrawfordWilliams1985}. Furthermore, the normalised left eigenvector components are the reciprocals of the normalised right eigenvector components \cite[Theorem~7.33]{Saaty1980}.
The two methods are also equivalent for consistent matrices.

However, they may differ if the number of alternatives is at least four and the matrix is inconsistent.
To the best of our knowledge, this problem was first presented by \cite{JohnsonBeineWang1979}.

\begin{example} \label{Examp1}
\cite{JohnsonBeineWang1979}
Consider the following pairwise comparison matrix with four alternatives:
\[
\mathbf{A} = \left[
\begin{array}{K{2em} K{2em} K{2em} K{2em}}
    1     	& 3		  	& 1/3   	& 1/2   \\
    1/3		& 1       	& 1/6		& 2 \\
    3		& 6			& 1      	& 1 \\
    2	 	& 1/2	  	& 1			& 1 \\
\end{array}
\right].
\]
Its right eigenvector, with the sum of weights normalised to 100, is:
\[
\mathbf{w}^{(R)} \left( \mathbf{A} \right) = \left[
\begin{array}{K{3em} K{3em} K{3em} K{3em}}
    18.44 & 15.19 & 43.64 & 22.73 \\
\end{array}
\right].
\]
Hence, the fourth alternative is ranked above the first.

The left eigenvector equals
\[
\mathbf{w}^{(L)} \left( \mathbf{A} \right) = \left[
\begin{array}{K{3em} K{3em} K{3em} K{3em}}  
    24.82 & 38.78 & 10.49 & 25.91 \\
\end{array}
\right],
\]
and its entry-wise reciprocal is
\[
\mathbf{w}^{(-L)} \left( \mathbf{A} \right) = \left[
\begin{array}{K{3em} K{3em} K{3em} K{3em}}
    20.14 & 12.89 & 47.67 & 19.29 \\
\end{array}
\right].
\]
Clearly, the weights from $\mathbf{w}^{(R)} $ and $\mathbf{w}^{(-L)} $ are different (for instance, the priority of the second alternative is higher by more than 15\% according to the right eigenvector), which affects even their ordering: in contrast to the right eigenvector, the inverse left eigenvector prefers the first alternative to the fourth.
\end{example}

In Example~\ref{Examp1}, the consistency ratio of matrix $\mathbf{A}$ is $\mathit{CR} \left( \mathbf{A} \right) \approx 0.331$.
But the right and inverse left eigenvectors may lead to a different ranking even if the inconsistency ratio is acceptable according to the recommendation of Saaty ($\mathit{CR}(\mathbf{A}) < 0.1$).

\begin{example} \label{Examp2}
\cite{DoddDoneganMcMaster1995}
Consider the following pairwise comparison matrix with five alternatives:
\[
\mathbf{B} = \left[
\begin{array}{K{2em} K{2em} K{2em} K{2em} K{2em}}
    1     & 1     & 3     & 9     & 9     \\
    1     & 1     & 5     & 8     & 5     \\
     1/3  &  1/5  & 1     & 9     & 5     \\
     1/9  &  1/8  &  1/9  & 1     & 1     \\
     1/9  &  1/5  &  1/5  & 1     & 1     \\
\end{array}
\right],
\]
where $\mathit{CR} \left( \mathbf{B} \right) \approx 0.082$.
Its right eigenvector is
\[
\mathbf{w}^R \left( \mathbf{B} \right) = \left[
\begin{array}{K{3em} K{3em} K{3em} K{3em} K{3em}}
    36.57 & 38.96 & 16.72 & 3.47 & 4.29 \\
\end{array}
\right],
\]
while the entry-wise reciprocal of the left eigenvector is
\[
\mathbf{w}^{(-L)} \left( \mathbf{B} \right) = \left[
\begin{array}{K{3em} K{3em} K{3em} K{3em} K{3em}}
    40.64 & 36.42 & 15.07 & 3.44 & 4.43 \\
\end{array}
\right].
\]
Consequently, a rank reversal emerges with respect to the most important alternative, which is either the first (inverse left eigenvector), or the second (right eigenvector).
\end{example}

Examples~\ref{Examp1} and \ref{Examp2} imply the question: Does the set of $n \times n$ pairwise comparison matrices whose left and right eigenvectors are reciprocals coincide with the set of consistent matrices for $n \geq 4$?
The answer is negative.

\begin{example} \label{Examp3}
\cite{DeTurck1987}
Consider the following pairwise comparison matrix with four alternatives:
\[
\mathbf{C} = \left[
\begin{array}{K{2em} K{2em} K{2em} K{2em}}
    1     &   8/5 &  1/4  & 4     \\
     5/8  & 1     &  5/8  & 10    \\
    4     &   8/5 & 1     & 4     \\
     1/4  & 1/10  &  1/4  & 1     \\
\end{array}
\right].
\]
The right eigenvector equals
\[
\mathbf{w}^R \left( \mathbf{C} \right) = \left[
\begin{array}{K{3em} K{3em} K{3em} K{3em} K{3em}}
      2/9  &   5/18 &   4/9  &   1/18 \\
\end{array}
\right],
\]
while the left eigenvector equals
\[
\mathbf{w}^{(L)} \left( \mathbf{C} \right) = \left[
\begin{array}{K{3em} K{3em} K{3em} K{3em} K{3em}}
      1/4  &   1/5 &   1/8  &   1 \\
\end{array}
\right].
\]
It can be checked that $\mathbf{w}^{(-L)} \left( \mathbf{C} \right) = \mathbf{w}^R \left( \mathbf{C} \right)$, even though $c_{12} c_{23} = 1 \neq c_{13}$.
\end{example}

\cite{DeTurck1987} also constructs, for larger $n$, inconsistent matrices---in which all entries are ones except for one entry $\alpha$ and another entry $1 / \alpha$ in each column and row---for which the two approaches give the same weight vector.

Three studies investigate the differences between the right and inverse left eigenvectors via Monte Carlo simulations.

\cite{IshizakaLusti2006} generates five groups of inconsistent matrices for each $3 \leq n \leq 7$: $0.02(k-1) \leq \mathit{CR} < 0.02k$ holds in group $k$, $1 \leq k \leq 5$. In particular, the $n-1$ entries directly above the main diagonal are randomly selected from the Saaty scale~\eqref{Saaty_scale}. All other entries are derived from these values to get a consistent matrix, but then some entries are perturbed as follows. The number of perturbed entries is uniformly chosen from the set of integers between 0 (no entry is perturbed) and $n(n-1)/2$ (the number of entries in the upper triangle). The maximal shift of any entry is four scale positions. The probability of each new value is 1/9, and the original comparison receives the remaining probability.
For instance, if the entry to be perturbed is 7, the new comparison can be 3, 4, 5, 6, 8, 9, each occurring with a chance of 1/9, while it remains 7 with a probability of 1/3. The choice of the entry to be perturbed is random; an entry may be changed more than once. Finally, when the required number of perturbations is carried out, the inconsistency ratio $\mathit{CR}$ is computed, and the matrix is either assigned to the corresponding group $k$ or dismissed if $\mathit{CR} \geq 0.1$. The priorities are calculated with both the right and inverse left eigenvectors, as well as with the row geometric mean. One hundred matrices are considered for each value of $n$ and group $k$, resulting in 2500 matrices in total.

Based on this simulation experiment, \cite{IshizakaLusti2006} concludes that
(1) the differences between the rankings of different weighting methods are slight;
(2) only close priorities suffer from ranking contradictions;
(3) the probability of a ranking contradiction increases with the level of inconsistency and the number of alternatives.

\cite{BozokiRapcsak2008} generates 100 million pairwise comparison matrices for $n=5$ to detect rank reversals in the weights computed from the right and inverse left eigenvectors. Each entry above the diagonal is uniformly and randomly chosen from the Saaty scale~\eqref{Saaty_scale}, analogous to the definition of the random index $\mathit{RI}$. The frequency of rank reversal is found to increase as a function of the inconsistency ratio, although the case $\mathit{CR} > 0.19$ is not analysed.

Currently, \cite{Csato2024a} provides probably the most extensive numerical study on right-left asymmetry. The paper adopts the method of \cite{SzadoczkiBozokiJuhaszKadenkoTsyganok2023} to obtain pairwise comparison matrices:
\begin{enumerate}
\item
A uniformly distributed random number $w_i$ is chosen from the interval $\left[ 1,\, 9 \right]$ for each $1 \leq i \leq n$.
\item
The consistent pairwise comparison matrix $\mathbf{A} = \left[ a_{ij} = w_i / w_j \right]$ is computed.
\item
For all $i \neq j$, either $a_{ij}$ or $a_{ji} = 1/a_{ij}$ is perturbed depending on which entry is greater. \\
If $a_{ij} \geq 1$, the perturbed entry $\tilde{a}_{ij}$ equals
\begin{equation} \label{Eq_perturbation}
\tilde{a}_{ij} =
\begin{cases}
    a_{ij} + \varepsilon_{ij} & \text{if } a_{ij} + \varepsilon_{ij} \geq 1 \\
    1 / \left[ 1 - \varepsilon_{ij} - \left( a_{ij} - 1 \right) \right] & \text{otherwise},
\end{cases}
\end{equation}
where $\varepsilon_{ij}$ is a uniformly distributed random value from the interval $\left[ -\Delta;\, \Delta \right]$. \\
If $a_{ij} < 1$, then $a_{ji} > 1$, and $\tilde{a}_{ji}$ is computed according to \eqref{Eq_perturbation}.
\item
The reciprocity of the matrix is retained by modifying the pair of the perturbed entry.
\end{enumerate}
The main advantage of this method resides in the uniform distribution of the perturbed entries on the Saaty scale~\eqref{Saaty_scale}, which is in line with the approach of \cite{IshizakaLusti2006}, but the latter only allows for integers from 1 to 9 and their reciprocals.

\cite{Csato2024a} considers three values (1, 2, 3) of parameter $\Delta$, while the number of alternatives is $4 \leq n \leq 9$. In each case, one million matrices are generated, giving a full sample of 18 million matrices. For any matrix, the inconsistency ratio $\mathit{CR}$ and the weight vectors $\mathbf{w}^R$, $\mathbf{w}^{(-L)}$, and $\mathbf{w}^{(\mathit{GM})}$ are computed such that the weights sum up to one. These vectors are compared by four different metrics:
(a) Euclidean distance (the root of the sum of squared differences);
(b) Chebyshev distance (the absolute value of the maximal difference);
(c) maximal ratio (the maximal ratio of the weights for the same alternative);
(d) Kendall tau (the similarity of the orderings derived from the weight vectors).

Row geometric mean turns out to be a good compromise between the right and inverse left eigenvectors, as it is close to the vector given by multiplying the right and inverse left eigenvectors componentwise. This also remains true from an ordinal perspective: the ranking derived from the row geometric mean is not farther from the ranking derived from the right eigenvector than the ranking derived from the inverse left eigenvector with a probability of at least 95\% (except for four alternatives when this result holds only if $\mathit{CR} < 0.2$).
In contrast to the finding of \cite{IshizakaLusti2006}, the differences between the three weighting methods do not always increase with the number of alternatives $n$: the mean Kendall tau is about the same in a given inconsistency interval if there are at least seven alternatives.

Finally, \cite{Csato2024a} presents three examples that might inspire further studies.
First, rank reversal between the right and inverse left eigenvectors may emerge even for a slightly inconsistent matrix.

\begin{example} \label{Examp4}
\cite{Csato2024a}
Consider the following pairwise comparison matrix with four alternatives:
\[
\mathbf{D} = \left[
\begin{array}{K{4em} K{4em} K{4em} K{4em}}
    1     & 0.4759 & 0.9832 & 0.4025 \\
    2.1011 & 1     & 1.9975 & 0.7374 \\
    1.0171 & 0.5006 & 1     & 0.3704 \\
    2.4842 & 1.3560 & 2.6998 & 1 \\
\end{array}
\right].
\]
The right eigenvector is
\[
\mathbf{w}^R \left( \mathbf{D} \right) = \left[
\begin{array}{K{4em} K{4em} K{4em} K{4em}}
    15.042 & 30.274 & 15.037 & 39.647 \\
\end{array}
\right],
\]
and the inverse left eigenvector is
\[
\mathbf{w}^{(-L)} \left( \mathbf{D} \right) = \left[
\begin{array}{K{4em} K{4em} K{4em} K{4em}}
    15.036 & 30.281 & 15.049 & 39.635 \\
\end{array}
\right].
\]
Thus, a rank reversal arises with respect to the first and third alternatives, although $\mathit{CR} \left( \mathbf{D} \right) \approx 0.0007$.
\end{example}

According to Example~\ref{Examp4}, a class of pairwise comparison matrices with arbitrarily small inconsistency may exist, for which the right and inverse left eigenvectors result in a different ranking of the alternatives.

Second, the right and inverse left eigenvectors could imply a fully reversed order of the alternatives.

\begin{example} \label{Examp5}
\cite{Csato2024a}
Consider the following pairwise comparison matrix with five alternatives:
\[
\mathbf{E} = \left[
\begin{array}{K{4em} K{4em} K{4em} K{4em} K{4em}}
    1     & 1.624 & 0.574 & 1.072 & 1.054 \\
    0.616 & 1     & 1.132 & 1.089 & 1.269 \\
    1.743 & 0.884 & 1     & 1.515 & 0.467 \\
    0.933 & 0.919 & 0.660 & 1     & 1.694 \\
    0.949 & 0.788 & 2.140 & 0.590 & 1 \\
\end{array}
\right].
\]
The right eigenvector is
\[
\mathbf{w}^R \left( \mathbf{E} \right) = \left[
\begin{array}{K{4em} K{4em} K{4em} K{4em} K{4em}}
    19.75 & 19.16 & 20.85 & 19.53 & 20.71 \\
\end{array}
\right],
\]
and the inverse left eigenvector is
\[
\mathbf{w}^{(-L)} \left( \mathbf{E} \right) = \left[
\begin{array}{K{4em} K{4em} K{4em} K{4em} K{4em}}
    20.25 & 20.55 & 19.31 & 20.27 & 19.62 \\
\end{array}
\right].
\]
The right eigenvector gives the ranking $3 \succ 5 \succ 1 \succ 4 \succ 2$ among the alternatives, which is reversed to $3 \prec 5 \prec 1 \prec 4 \prec 2$ based on the inverse left eigenvector.
\end{example}

It remains an open question whether an analogous example exists for four alternatives. In addition, Example~\ref{Examp5} calls for finding a class of pairwise comparison matrices for which the right and inverse left eigenvectors lead to a reversed order of the alternatives.

\cite[p.~398]{IshizakaLusti2006} states that ``\emph{only very close priorities suffer from ranking contradictions}''. However, the following example suggests otherwise.

\begin{example} \label{Examp6}
\cite{Csato2024a}
Consider the following pairwise comparison matrix with five alternatives:
\[
\mathbf{F} = \left[
\begin{array}{K{4em} K{4em} K{4em} K{4em} K{4em}}
    1     & 0.371 & 2.013 & 5.389 & 0.243 \\
    2.698 & 1     & 4.596 & 7.527 & 0.736 \\
    0.497 & 0.218 & 1     & 2.321 & 0.167 \\
    0.186 & 0.133 & 0.431 & 1     & 0.385 \\
    4.120 & 1.359 & 5.973 & 2.598 & 1 \\
\end{array}
\right].
\]
The right eigenvector is
\[
\mathbf{w}^R \left( \mathbf{F} \right) = \left[
\begin{array}{K{4em} K{4em} K{4em} K{4em} K{4em}}
    15.26 & 33.23 & 7.74  & 5.68  & 38.08 \\
\end{array}
\right],
\]
and the inverse left eigenvector is
\[
\mathbf{w}^{(-L)} \left( \mathbf{F} \right) = \left[
\begin{array}{K{4em} K{4em} K{4em} K{4em} K{4em}}
    15.29 & 37.84 & 8.55  & 4.93  & 33.39 \\
\end{array}
\right].
\]
The absolute value of the difference between the weights of the second (the best based on the inverse left eigenvector) and the fifth (the best based on the right eigenvector) alternatives is $4.85$ and $4.44$, respectively, when the sum of weights is 100, which is substantial relative to the weights.
Meanwhile, the inconsistency ratio $\mathit{CR}$ remains below 0.1.
\end{example}

\section{Group incoherence for choice} \label{Sec4}

In group decision-making, the alternatives have to be ranked based on the preferences of more than one decision-maker, expressed by pairwise comparison matrices. This requires aggregating either the individual weights or the individual pairwise comparison matrices in order to get an overall assessment. In the latter case, the natural candidate is the (weighted) geometric mean of each entry: this is the unique quasi-arithmetic mean that satisfies reciprocity (hence, the aggregated matrix remains a pairwise comparison matrix) and positive homogeneity (multiplying individual preferences by the same positive scalar changes the aggregated preferences accordingly) \cite{AczelSaaty1983}. 

As a result, two different approaches can be considered to derive priorities in a group decision-making problem.
First, derive the priorities from the individual pairwise comparison matrices, and aggregate them. Second, aggregate the individual pairwise comparison matrices, and derive the priorities from this aggregated matrix.
\cite{PerezMokotoff2016} says that a weighting method satisfies \emph{group coherence for choice} if these two approaches always give the highest priority to the same alternative.

The row geometric mean method meets group coherence for choice \cite{Barzilai1997, PerezMokotoff2016}, but the eigenvector method violates it.

\begin{example} \label{Examp7}
\cite[Theorem~3.2]{PerezMokotoff2016}
Consider the following pairwise comparison matrices with five alternatives:
\[
\mathbf{A} = \left[
\begin{array}{K{3em} K{3em} K{3em} K{3em} K{3em}}
    1     & 1	  & 2	  & 1	  & 1 \\
    1	  & 1     & 1	  & 2	  & 1 \\
    1/2   & 1	  & 1     & 1	  & 2 \\
    1	  & 1/2	  & 1	  & 1     & 1 \\
    1	  & 1	  & 1/2	  & 1	  & 1 \\
\end{array}
\right]
\text{, and}
\]
\[
\mathbf{B} = \left[
\begin{array}{K{3em} K{3em} K{3em} K{3em} K{3em}}
    1     & 1	  & 1	  & 1	  & 2 \\
    1	  & 1     & 1	  & 2	  & 1 \\
    1	  & 1	  & 1     & 1	  & 1/2 \\
    1	  & 1/2	  & 1	  & 1     & 1 \\
    1/2	  & 1	  & 2	  & 1	  & 1 \\
\end{array}
\right].
\]
The right eigenvectors of them are
\[
\mathbf{w}^{(\mathit{EM})} \left( \mathbf{A} \right) = \left[
\begin{array}{K{3em} K{3em} K{3em} K{3em} K{3em}}
    23.06 & 22.43 & 20.26  & 17.02  & 17.23 \\
\end{array}
\right]
\text{, and}
\]
\[
\mathbf{w}^{(\mathit{EM})} \left( \mathbf{B} \right) = \left[
\begin{array}{K{3em} K{3em} K{3em} K{3em} K{3em}}
    23.06 & 22.43 & 17.23  & 17.02  & 20.26 \\
\end{array}
\right].
\]
Consequently, the first alternative is the best in both matrices $\mathbf{A}$ and $\mathbf{B}$.
On the other hand, if the matrices are aggregated by the entry-wise geometric mean method, then the right eigenvector of the aggregated matrix $\mathbf{C}$ equals
\[
\mathbf{w}^{(\mathit{EM})} \left( \mathbf{C} \right) = \left[
\begin{array}{K{3em} K{3em} K{3em} K{3em} K{3em}}
    22.69 & 23.12 & 18.38  & 17.42  & 18.38 \\
\end{array}
\right],
\]
that is, the second alternative receives the highest priority.

Note that the row geometric mean implies the ranking $1 \sim 2 \succ 3 \succ 4 \sim 5$ for matrix $\mathbf{A}$, $1 \sim 2 \succ 5 \succ 3 \sim 4$ for matrix $\mathbf{B}$, and $1 \sim 2 \succ 3 \sim 5 \succ 4$ for matrix $\mathbf{C}$, respectively.
\end{example}

\cite{Csato2017c} introduces the axiom \emph{inversion}. A weighting method satisfies inversion if the ranking of the alternatives derived from the transposed pairwise comparison matrix, which reverses all preferences of the decision-maker, is exactly the opposite of the original ranking. 

If a weighting method violates inversion, but satisfies the minimal requirements of anonymity (independence of the labelling of the alternatives), invariance to row multiplication (the relative weight of an alternative is multiplied by $\alpha$ if each entry in the corresponding row of the pairwise comparison matrix is multiplied by $\alpha$), then it also violates group coherence for choice \cite[Lemma~4.1 and Proposition~4.6]{Csato2017c}.

The eigenvector method violates inversion since the right eigenvector of the transposed pairwise comparison matrix is the inverse left eigenvector of the original matrix, as discussed in Section~\ref{Sec3}. Therefore, group incoherence for choice is closely related to right-left asymmetry.

The existing literature has not yet examined how likely it is that the eigenvector method does not satisfy group coherence for choice.

\section{Sensitivity to the intensity of preferences} \label{Sec5}

Let $\mathbf{A} = \left[ a_{ij} \right]$ be a pairwise comparison matrix and $\alpha > 0$ be a positive scalar.
Define the following set of pairwise comparison matrices:
\[
\mathcal{A} \left( \mathbf{A}, \alpha \right) = \left\{ \mathbf{B} = \left[ b_{ij} \right]: b_{ij} = a_{ij}^\alpha \text{ for all } 1 \leq i,j \leq n, \alpha > 0 \right\}.
\]
This class contains all pairwise comparison matrices derived from $\mathbf{A}$ by raising its entries to the positive power of $\alpha$.
For example, if the pairwise comparisons are only wins (represented by an arbitrary value $p$), draws (represented by 1), and losses (represented by $1/p$), it is far from clear how the value of $p$ should be chosen. Alternatively, the preferences may be given on a verbal scale, when their translation into numbers remains challenging \cite{CavalloIshizaka2023, SzadoczkiBozokiSiposGalambosi2025a}.

\cite{PetroczyCsato2021} introduces the axiom \emph{scale invariance}.
A weighting method satisfies scale invariance if, given any pairwise comparison matrix $\mathbf{A}$, every matrix in the set $\mathcal{A} \left( \mathbf{A}, \alpha \right)$ implies the same ranking of the alternatives. In other words, the ranking is not allowed to change if a different scale is used for pairwise comparisons. The row geometric mean method satisfies scale invariance \cite[Lemma~7]{PetroczyCsato2021}.

However, the eigenvector method violates scale invariance.

\begin{example} \label{Examp8}
\cite[Example~2.1]{GenestLapointeDrury1993}
Consider the following parametric pairwise comparison matrix with six alternatives:
\[
\mathbf{A}(p) = \left[
\begin{array}{K{3em} K{3em} K{3em} K{3em} K{3em} K{3em}}
    1     & p     & p     & 1/p   & p     & p \\
    1/p   & 1     & 1/p   & p     & p     & 1/p \\
    1/p   & p     & 1     & p     & p     & p \\
    p     & 1/p   & 1/p   & 1     & 1/p   & 1/p \\
    1/p   & 1/p   & 1/p   & p     & 1     & p \\
    1/p   & p     & 1/p   & p     & 1/p   & 1 \\
\end{array}
\right].
\]
This matrix, representing the results of a single round-robin chess tournament, is ordinally intransitive. For instance, Player 2 defeated Player 5, who won against Player 6, but Player 2 lost to Player 6. The matrix has five intransitive triads. Since the number of intransitive triads is at most $n(n^2 - 4)/24$ if $n$ is even \cite{KendallSmith1940}, which gives 8 for $n=6$, a considerable amount of intransitivity emerges in this competition.

The ranking by the row geometric mean method is $1 \sim 3 \succ 2 \sim 5 \sim 6 \succ 4$, independent of the value of parameter $p$. The top two players have four wins and one loss each, the middle three players have two wins and three losses each, while the weakest player has one win and four losses.

On the other hand, the ranking based on the right eigenvector is affected by parameter $p$. Players 2, 5, 6 always get the same weight because they all won against Player 4, lost to Players 1 and 3, and they have a circular triad among themselves.
If $p$ is smaller than a threshold around 3.5, Player 4 is ranked below these three players.
However, their order is reversed for a higher value of $p$, when the upset win of Player 4 against the strong Player 1 outweighs the worse overall record of Player 4.
\end{example}

An example with multiple rank reversals also exists.

\begin{example} \label{Examp9}
\cite[Example~6.2]{GenestLapointeDrury1993}
Consider the following parametric pairwise comparison matrix with seven alternatives:
\[
\mathbf{A}(p) = \left[
\begin{array}{K{3em} K{3em} K{3em} K{3em} K{3em} K{3em} K{3em}}
    1     & 1/p   & p     & p     & p     & p     & p \\
    p     & 1     & p     & p     & 1/p   & p     & 1/p \\
    1/p   & 1/p   & 1     & p     & p     & p     & p \\
    1/p   & 1/p   & 1/p   & 1     & p     & p     & p \\
    1/p   & p     & 1/p   & 1/p   & 1     & 1/p   & p \\
    1/p   & 1/p   & 1/p   & 1/p   & p     & 1     & p \\
    1/p   & 1/p   & 1/p   & 1/p   & 1/p   & 1/p   & 1 \\
\end{array}
\right].
\]

The ranking based on the row geometric mean is $1 \succ 2 \sim 3 \succ 4 \succ 5 \sim 6 \succ 7$, independent of parameter $p$ again.
If $p=2$, the inconsistency ratio is $\mathit{CR} \left( \mathbf{A}(p) \right) \approx 0.13$, and the ranking coincides with those from the row geometric mean. However, if the preferences are more intense with $p=4$, inconsistency becomes substantially higher ($\mathit{CR} \left( \mathbf{A}(p) \right) \approx 0.65$), and the ranking is $2 \succ 1 \succ 3 \succ 5 \succ 4 \succ 7 \succ 6$. Consequently, all ties are resolved, and the second alternative overtakes the first, the fifth overtakes the fourth, while the seventh overtakes the sixth, indicating a triple rank reversal.
\end{example}

Last but not least, scale invariance can be violated by the eigenvector method in real-world applications of pairwise comparison matrices.
\cite{PetroczyCsato2021} uses this methodology to allocate Formula One prize money among  the teams.
Team $i$ is said to score $k_{h,ij}$ points against team $j$ in race $h$ if there are $k$ instances in which a car of team $i$ finishes ahead of a car of team $j$. Since each team has (at most) two cars in every race, the maximum of $k$ in each race is four. The pairwise comparison matrix $\mathbf{A}^{(\alpha)}$ is defined by
\[
a_{ij}^{(\alpha)} = \left( \frac{ \sum_h k_{h,ij}}{\sum_{h} k_{h,ji}} \right)^{\alpha},
\]
where parameter $\alpha$ controls the inequality of the allocation.

The ranking of Formula One teams by their prize money remains unchanged as a function of $\alpha$ if the row geometric mean method is used. In the 2014 season, McLaren and Williams always receive the same amount, while Ferrari gets more than McLaren, and Marussia gets more than Sauber.
However, with the eigenvector method, McLaren receives less (more) than Williams if the value of $\alpha$ is below (above) 1.5. Furthermore, for $\alpha$ higher than 0.11, McLaren is preferred to Ferrari. The difference between the shares of Marussia and Sauber is even non-monotonic \cite[Figure~2]{PetroczyCsato2021}: it is first increasing and positive, then decreasing and becomes negative, and, finally, again increasing but remains negative when $\alpha$ grows.

\section{Non-monotonicity} \label{Sec6}

Entry $a_{ij}$ of a pairwise comparison matrix $\mathbf{A} = \left[ a_{ij} \right]$ quantifies the dominance of alternative $i$ over alternative $j$; hence, increasing $a_{ij}$ is expected to be beneficial for alternative $i$. Inspired by this idea, \cite{CsatoPetroczy2021} considers two monotonicity properties.

Let $\mathbf{A} = \left[ a_{ij} \right]$ be any pairwise comparison matrix, $i \neq j$ two different alternatives, and $\mathbf{A}^{\prime} = \left[ a_{ij}^{\prime} \right]$ be an arbitrary pairwise comparison matrix identical to $\mathbf{A}$ except for $a_{ij}^{\prime} > a_{ij}$ and, naturally, $a_{ji}^{\prime} < a_{ji}$ because of reciprocity.
A weighting method is said to be \emph{rank monotonic} if $\mathbf{w} \left( \mathbf{A} \right)_i \geq \mathbf{w} \left( \mathbf{A} \right)_k$ implies $\mathbf{w} \left( \mathbf{A}^{\prime} \right)_i \geq \mathbf{w} \left( \mathbf{A}^{\prime} \right)_k$ for all $1 \leq k \leq n$.

Rank monotonicity requires that, if an alternative is originally ranked at least as high as another alternative, and one pairwise comparison changes in favour of it, then its rank does not worsen with respect to the other alternative. The motivation is clear: any counter-intuitive rank reversal would be against the intention of the decision-maker.

Rank monotonicity, which focuses on the ranking of the alternatives, has an analogous form regarding the (normalised) weights.
A weighting method is called \emph{weight monotonic} if the inequality
\[
\frac{\mathbf{w} \left( \mathbf{A}^{\prime} \right)_i}{\sum_{k=1}^n \mathbf{w} \left( \mathbf{A}^{\prime} \right)_k} \geq \frac{\mathbf{w} \left( \mathbf{A} \right)_i}{\sum_{k=1}^n \mathbf{w} \left( \mathbf{A} \right)_k}
\]
holds for any pairwise comparison matrices $\mathbf{A}$, $\mathbf{A}^{\prime}$, and alternative $i$.
In other words, increasing $a_{ij}$ is not allowed to decrease the normalised weight of alternative $i$.

The row geometric mean method satisfies both rank monotonicity and weight monotonicity \cite[Proposition~1]{CsatoPetroczy2021}: when $a_{ij}$ becomes higher, the product of the $i$th row increases, the product of the $j$th row decreases, while the product of any other rows does not change.

On the other hand, the eigenvector method does not satisfy these monotonicity axioms.

\begin{example} \label{Examp10}
\cite[Example~1]{CsatoPetroczy2021}
Take the following parametric pairwise comparison matrix of order six:
\[
\mathbf{A}(\beta) = \left[
\begin{array}{K{3em} K{3em} K{3em} K{3em} K{3em} K{3em}}
    1     & \beta     & 8     &  1/9  &  1/4  &  1/9 \\
    1/ \beta   & 1     &  1/5  &  1/9  & 1     &  1/4 \\
     1/8  & 5     & 1     &  1/8  &  1/2  &  1/7 \\
    9     & 9     & 8     & 1     & 7     & 8     \\
    4     & 1     & 2     &  1/7  & 1     &  1/9 \\
    9     & 4     & 7     &  1/8  & 9     & 1     \\
\end{array}
\right].
\]
$\mathbf{w}^{(\mathit{EM})} \left( \mathbf{A}(0.3) \right)_1 > \mathbf{w}^{(\mathit{EM})} \left( \mathbf{A}(0.3) \right)_5$, but $\mathbf{w}^{(\mathit{EM})} \left( \mathbf{A}(0.5) \right)_1 < \mathbf{w}^{(\mathit{EM})} \left( \mathbf{A}(0.5) \right)_5$, which shows the violation of rank monotonicity when the value of $\beta$ increases from 0.3 to 0.5.
In addition, the normalised weight of the first alternative has a local minimum around $\alpha \approx 0.5$ using the eigenvector method; hence, weight monotonicity does not hold.
\end{example}

\cite{CsatoPetroczy2021} has also conducted a simulation analysis to explore how likely the non-monotonicity of the eigenvector method is. Following the methodology of \cite{BozokiRapcsak2008}, a large number of pairwise comparison matrices have been generated: all entries in the upper triangle are chosen independently and uniformly from the Saaty scale~\eqref{Saaty_scale}. For each pairwise comparison matrix, $n(n-1) / 2$ perturbed matrices have been considered by substituting one entry above the diagonal by the next element from the Saaty scale (and by 10 if it is 9). Finally, the probability that rank or weight monotonicity is violated for a given matrix with an inconsistency ratio $\mathit{CR}$ between $0.01(m-1)$ and $0.01m$ has been computed for all values of $m$.

Obviously, the eigenvector method satisfies both monotonicity properties if the number of alternatives is $n=3$, because it coincides with the row geometric method then \cite[p.~393]{CrawfordWilliams1985}.
For $n=4$, all of the $17^6 = 24{,}137{,}569$ matrices have been checked according to the procedure above, but no violation of either rank or weight monotonicity has been found \cite{CsatoPetroczy2021}.
10 million random matrices have been generated for $5 \leq n \leq 8$, and 5 million for $n=9$. The frequency of violating rank monotonicity increases roughly linearly with the inconsistency of the pairwise comparison matrix. No example has appeared if $\mathit{CR} < 0.2$, and the problem emerges with a probability of around 2\% for heavily inconsistent matrices.
Weight monotonicity is violated the most frequently around an inconsistency ratio of 0.5. In some cases, both monotonicity properties are violated, similar to Example~\ref{Examp10}.

Nevertheless, a number of open questions remain about the non-monotonicity of the eigenvector method.
\cite{CsatoPetroczy2021} does not provide any example if $n=4$, or if the inconsistency ratio remains below 0.1. The connection between non-monotonicity and right-left asymmetry has not been investigated yet. Non-monotonicity can be quantified not only by the probability of matrices for which the required implication does not hold, but also by the extent to which the rank and the normalised weight of alternative $i$ can worsen following a given increase of entry $a_{ij}$.

\section{Pareto (in)efficiency} \label{Sec7}

Let $\mathbf{A} = \left[ a_{ij} \right]$ be a pairwise comparison matrix.  
A weight vector $\mathbf{w} = \left[ w_i \right]$ is called \emph{efficient} for the multi-objective optimization problem
\[
\min \limits_{\mathbf{x} \in \mathbb{R}^n_{++}}
\left|  a_{ij}  - \frac{x_i}{x_j} \right|_{1 \leq i \neq j \leq n}
\]
if there does not exist any positive weight vector $\mathbf{w}^{\prime} = \left[ w_{i}^{\prime} \right]$ that Pareto improves the approximation, namely, it is at least as good in approximating all pairwise comparisons, but strictly better for at least one:
\begin{align} 
 \left| a_{ij} - \frac{w^{\prime}_i}{w^{\prime}_j} \right| &\leq \left| a_{ij} - \frac{w_i}{w_j} \right| \qquad \text{ for all } 1 \leq i,j \leq n, \label{Eff_prop1} \\   
 \left| a_{k{\ell}} - \frac{w^{\prime}_k}{w^{\prime}_{\ell}} \right| & < \left| a_{k{\ell}} - \frac{w_k}{w_{\ell}} \right| \qquad \text{ for some } 1 \leq k,\ell \leq n. \label{Eff_prop2}   
\end{align}
A weight vector is called \emph{inefficient} if it does not provide an efficient approximation.

If a weight vector $\mathbf{w}$ is inefficient, and another weight vector $\mathbf{w^{\prime}}$ satisfies \eqref{Eff_prop1}--\eqref{Eff_prop2},
$\mathbf{w^{\prime}}$ is said to \emph{dominate} $\mathbf{w}$.

The right eigenvector can be \emph{inefficient}.

\begin{example} \label{Examp11}
Consider the following pairwise comparison matrix:
\[
\mathbf{A} = \left[
\begin{array}{K{3em} K{3em} K{3em} K{3em}}
    1     & 2     & 4     &  7    \\
    1/2   & 1     & 3     &  2    \\
    1/4   & 1/3   & 1     &  3    \\
    1/7   & 1/2   & 1/3   &  1    \\
\end{array}
\right].
\]
Its inconsistency ratio is $\mathit{CR}(\mathbf{A}) = 0.0741$, and the right eigenvector is
\[
\mathbf{w}^{(\mathit{EM})} = \left[
\begin{array}{K{5em}}
\textcolor{red}{51.5592} \\
26.3375 \\
14.1984 \\
\textcolor{white}{0}7.9049
\end{array}
\right].
\]
This gives the following approximation:
\[
\left[ \frac{{w}^{(\mathit{EM})}_i}{{w}^{(\mathit{EM})}_j} \right] = \left[
\begin{array}{K{4em} K{4em} K{4em} K{4em}}
1 & \textcolor{red}{1.9576} & \textcolor{red}{3.6313} & \textcolor{red}{6.5224} \\
\textcolor{red}{0.5108} & 1 & 1.8550 & 3.3318 \\
\textcolor{red}{0.2754} & 0.5391 & 1 & 1.7962 \\
\textcolor{red}{0.1533} & 0.3001 & 0.5567 & 1 \\
\end{array}
\right].
\]

However, the weight vector
\[
\mathbf{w}^{\prime} = \left[
\begin{array}{K{5em}}
52.0937\\
26.0469\\
14.0418\\
\textcolor{white}{0}7.8177
\end{array}
\right] = 
0.988966 \cdot \left[
\begin{array}{K{5em}}
\textcolor{ForestGreen}{52.6749} \\
26.3375 \\
14.1984 \\
\textcolor{white}{0}7.9049
\end{array}
\right]
\]
provides
\[
\left[ \frac{{w}^{\prime}_i}{{w}^{\prime}_j} \right] = \left[
\begin{array}{K{4em} K{4em} K{4em} K{4em}}
1 & \textcolor{blue}{2} & \textcolor{ForestGreen}{3.7099} & \textcolor{ForestGreen}{6.6636} \\
\textcolor{blue}{1/2} & 1 & 1.8550 & 3.3318 \\
\textcolor{ForestGreen}{0.2695} & 0.5391 & 1 & 1.7962 \\
\textcolor{ForestGreen}{0.1501} & 0.3001 & 0.5567 & 1 \\
\end{array}
\right].
\]
Consequently, inequality~\eqref{Eff_prop1} holds for all $1 \leq i,j \leq 4$, and the strict inequality~\eqref{Eff_prop2} holds for six pairwise comparisons when $(k,\ell) \in \{ (1,2), (1,3), (1,4), (2,1), (3,1), (4,1) \}$.
For instance, $\left| w^{\prime}_1 / w^{\prime}_2 - a_{12} \right| = \left| 2-2 \right| = 0 < \left| w^{(\mathit{EM})}_1 / w^{(\mathit{EM})}_2 - a_{12} \right| = \left| 2 - 1.9576 \right| = 0.0424$.
Weight vector $\mathbf{w}^{\prime}$ dominates $\mathbf{w}^{(\mathit{EM})}$.
\end{example}

Furthermore, the right eigenvector may be inefficient  even if the inconsistency of the pairwise comparison matrix is arbitrarily small.

\begin{example}  \label{Examp12}
\cite{Bozoki2014a}
Let $n \geq 4$ and $\mathbf{A}(p,q)$ be a pairwise comparison matrix of order $n$ such that
\[
\mathbf{A}(p,q) = \left[
\begin{array}{K{3em} K{3em} K{3em} K{3em} K{3em} K{3em} K{3em}}
   1     &     p    &     p    &   p     &  \ldots &    p    &  p   \\
 1/p     &     1    &     q    &   1     &  \ldots &    1    &  1/q     \\
 1/p     &    1/q   &     1    &   q     &  \ldots &    1    &  1       \\
  \vdots &   \vdots &   \vdots &  \ddots &         &  \vdots &   \vdots  \\
  \vdots &   \vdots &   \vdots &         &  \ddots &  \vdots &   \vdots  \\
 1/p     &     1    &     1    &   1     &  \ldots &    1    &   q       \\
 1/p     &     q    &     1    &   1     &  \ldots &    1/q    &   1
\end{array}
\right],
\]
where $p$ and $q$ are positive numbers. $\mathit{CR} \left( \mathbf{A}(p,q) \right)$ goes to zero as $q$ goes to 1, and the right eigenvector of $\mathbf{A}(p,q)$ remains inefficient \cite[Proposition~2.4]{Bozoki2014a}.
\end{example}

In some special cases, the efficiency of the right eigenvector is guaranteed.
The eigenvector is always Pareto efficient if the pairwise comparison matrix differs from a consistent one only in one entry (and its reciprocal) \cite{Abele-NagyBozoki2016, daCruzFernadesFurtado2021}. 
\cite{Abele-NagyBozokiRebak2018} and \cite{Furtado2023} extend this result to double-perturbed matrices that can be made consistent by changing two elements and their reciprocals.  
\cite{FernandesFurtado2022} finds a class of triple perturbed pairwise comparison matrices for which the principal right eigenvector is efficient: a $2 \times 2$ submatrix contains three perturbed entries and does not contain a diagonal entry, while the magnitude of the perturbations is constrained. \cite{FernandesPalheira2024} continues the analysis of these triple perturbed matrices and determines which magnitudes of perturbations characterise the efficiency of the right eigenvector.
However, the right eigenvector cannot be efficient for all triple perturbed pairwise comparison matrices, as any $4 \times 4$ matrix is triple perturbed, and Example~\ref{Examp11} shows such a matrix with a Pareto inefficient eigenvector.

Interestingly, the geometric mean of the right eigenvector $\mathbf{w}^{(R)}$ and the entry-wise reciprocal of the left eigenvector $\mathbf{w}^{(-L)}$ is always efficient \cite[Theorem~20]{FurtadoJohnson2026}.

\cite{FurtadoJohnson2026} also studies the empirical probability that the right eigenvector is efficient based on randomly generated pairwise comparison matrices. First, an $n \times n$ matrix $\mathbf{C} = \left[ c_{ij} \right]$ is created, where $c_{ij}$ is drawn uniformly from the unit interval. The pairwise comparison matrix $\mathbf{A} = \left[ a_{ij} \right]$ is defined by $a_{ij} = c_{ij} / c_{ji}$ to guarantee reciprocity. For each value of $4 \leq n \leq 25$, the authors count the number of matrices for which the eigenvector is efficient in their sample of 50 thousand random matrices.
In the case of four, five, and six alternatives, the probability of an inefficient right eigenvector turns out to be around 10\%. This proportion gradually decreases as parameter $n$ grows, and remains below 1\% if $n \geq 17$.

\cite{DulebaMoslem2019} examines the Pareto-optimality of the eigenvector for pairwise comparison matrices coming from an AHP survey in public transport. Even though the dataset contains only one matrix with more than three alternatives, its eigenvector is inefficient; thus, a dominating weight vector is computed. At first sight, this approach seems to solve the issue of Pareto inefficiency in applications---but the dominating vector is usually non-unique, and choosing among them is far from straightforward.



The row geometric mean weight vector $\mathbf{w}^{(\mathit{GM})}$ is always efficient \cite[Corollary~7]{BlanqueroCarrizosaConde2006}.
On the other hand, finding a necessary and sufficient condition for the efficiency of the eigenvector in the general case remains an open problem.

\section{Conclusions} \label{Sec8}

We have presented five potential shortcomings of the eigenvector method, the most popular weighting procedure for multiplicative pairwise comparison matrices. All of them are satisfied by the row geometric mean method, which gives a powerful argument for its use.

Nevertheless, violations of these properties imply only that the eigenvector method is a poor choice in certain cases. Even though there exist some simulation studies \cite{CsatoPetroczy2021, Csato2024a, FurtadoJohnson2026}, it remains largely unknown whether these problems emerge only for ``pathological'' matrices (e.g.\ with a high level of inconsistency), or whether they are frequent and troubling phenomena in practice. It would also be interesting to see how other weighting methods perform with respect to the theoretical requirements above.

Another promising direction for future studies is exploring how the violation of the last four properties (group coherence for choice, scale invariance, monotonicity, Pareto efficiency) is related to right-left asymmetry. Section~\ref{Sec4} has revealed that the difference between the right and inverse left eigenvectors is connected to group incoherence for choice, and the important result \cite[Theorem~20]{FurtadoJohnson2026} uncovers the efficiency of the geometric mean of the right and inverse left eigenvectors. Based on these findings, some theoretical weaknesses of the eigenvector method could be rooted in the asymmetric choice of the right rather than the inverse left eigenvector.

\begin{acknowledgement}
The research was supported by the National Research, Development and Innovation Office under Grants Advanced 152220 and FK 145838, and the J\'anos Bolyai Research Scholarship of the Hungarian Academy of Sciences.
\end{acknowledgement}
\ethics{Competing Interests}{The authors have no conflicts of interest to declare that are relevant to the content of this chapter.}

\eject

\bibliographystyle{spmpsci}
\bibliography{All_references}

\end{document}